\documentclass{elsart}
\usepackage{epsfig}

\newtheorem{theorem}{Theorem}[section]

\newtheorem{lemma}{Lemma}[section]

\usepackage{amssymb}
\usepackage{amsmath} 
\usepackage{epsfig}

\begin{document}

\begin{frontmatter}

\title{\bf  Neural networks  with non-smooth and impact activations}
\author{M. U. Akhmet$^{a}$\corauthref{3}} \ead{marat@metu.edu.tr} and \author{E. Y{\i}lmaz$^b$} \ead{enes@metu.edu.tr}

\address{$^a$Department of Mathematics, Middle East
Technical University, 06531 Ankara, Turkey}
\address{$^b$Institute of Applied Mathematics, Middle East
Technical University, 06531 Ankara, Turkey}
\corauth[3]{Corresponding author}

\begin{abstract}
In this paper,  we  consider  a model of impulsive recurrent neural networks  with piecewise constant delay.  The dynamics are presented by differential equations with discontinuities such as impulses at fixed moments and  piecewise constant argument of generalized type.  Sufficient conditions ensuring the existence, uniqueness and global asymptotic stability of  the equilibrium point  are obtained. By employing Green's function  we derive new result  of  existence  of the periodic solution. The  global asymptotic stability  of the solution is investigated. Examples with numerical simulations are given to validate the theoretical  results.
\end{abstract}
\begin{keyword} 
Recurrent neural network; Impulsive differential equation; Piecewise constant delay; Equilibrium; Periodic solution; Global asymptotic stability

%{\it 2000 MSC:} 
\end{keyword}

\end{frontmatter}

\section{Introduction}\label{secintro}
Recurrent neural networks  and impulsive recurrent neural networks  have been investigated due to their extensive applications in classification of patterns, associative memories, image processing, optimization problems, and other areas \cite{hop,CNNT,CNNA,CNNA1,mfp,cgp,coombes,HNNWI1,HNNWI3,HNNWI4,HNNWDI1,HNNWDI2,HNNWDI3,period}. It is well  known that these  applications depend crucially on the dynamical behavior of the networks. For example, if a neural network is employed to solve some optimization problems, it is highly desirable for the neural network
to have a unique globally stable equilibrium \cite{sabri,Cao4,Cao5,Cao6,xcl,hcw,pvd,zeng}. Therefore, stability analysis of neural networks has received much attention and various stability conditions have been obtained over the past years. Another interesting subject is to study the dynamical behavior of existence  of the periodic solutions in  recurrent neural networks. These periodic solutions present periodic pattern and have been used in learning theory, which are meant to capture the idea that certain activities or motions are learned by repetition \cite{town0,town}.

In numerical simulations and practical implementations of neural networks, it is essential to formulate a discrete-time system,   an analogue of the continuous-time system. Hence, stability for discrete-time  neural networks has also received
considerable attention from many researchers \cite{dtnns1,dtnns2,dtnns3,dtnns4,dtnns5,dtnns6}.
As we know, the reduction of differential equations with piecewise constant argument to discrete equations has  been  the main and possibly a unique way of stability analysis for these equations  \cite{wi3,wi4}.   Hence,  one has  not investigated the problem of stability completely, as only elements of a countable set were allowed to be discussed for initial moments. Finally,  only equations which  are  linear  with  respect to  the values of solutions at  non-deviated  moments of time   have been investigated.  That  narrowed significantly  the  class of systems. 
In papers \cite{ak0,ak2,ak3,ak4}, the theory  of differential equations with piecewise constant argument  has been generalized by Akhmet. Later, Akhmet gathered all  results for this type of differential equations in the book \cite{a0000}. All of these equations are reduced to equivalent integral equations such that one can investigate many problems, which have not been solved properly by using discrete equations, i.e., existence and uniqueness of solutions, stability and existence of periodic solutions. Moreover, since we do not need additional assumptions on the reduced discrete equations, the new method requires more easily verifiable conditions, similar to those for  ordinary differential equations.

In this paper, we develop the model of  recurrent neural networks to  differential equations with both impulses and  piecewise constant argument of generalized type. It is well known that impulsive differential equation \cite{ak00,SP,laks} is one of the basic instruments so the role of discontinuity has been understood better for the real world problems. In real world, many evolutionary processes are characterized by abrupt changes at certain time. These changes are called to be impulsive phenomena, which are included in many fields such as biology involving thresholds, bursting rhythm models,  physics, chemistry, population dynamics, models in economics, optimal control, etc.  In the literature, recurrent neural networks have been developed by implementing impulses and piecewise constant delay \cite{ak6,ak7,HNNWI1,HNNWI3,HNNWI4,HNNWDI1,HNNWDI2,HNNWDI3,period} issuing from different reasons: In implementation of electronic networks, the state of the networks is subject to instantaneous perturbations and experiences abrupt change at certain instants which may be caused by switching phenomenon, frequency change or other sudden noise. Furthermore, the dynamics of quasi-active dendrites with active spines is described by a system of point hot-spots (with an integrate-and-fire process), see \cite{dent1,dent2} for more details. This leads to the model of recurrent neural network with impulses. It is important to say that the neighbor moments of impulses may depend on each other. For example, the successive impulse moment may depend upon its predecessor. The reason for this phenomenon  is the interior design of a neural network. On the other hand, due to the finite switching speed of amplifiers and transmission of signals in electronic networks or finite speed for signal propagation in neural networks, time delays
exist \cite{CNNA,CNNA1,cgp,coombes}. Moreover,  the idea of involving delayed arguments in the recurrent neural networks can be explained by the fact  that we assume neural networks may {\it``memorize''} values of the phase variable at certain moments of time to utilize the values during middle process till the next moment. Thus, we arrive to differential equations with piecewise constant delay. Obviously,
the distances between the {\it``memorized''} moments may be very variative. Consequently, the concept of generalized type of piecewise constant argument is fruitful for recurrent neural networks \cite{a5,ak6,ak7}. 
Therefore, it is possible to apply  differential equations with both impulses and piecewise constant delay to neural networks theory.

The intrinsic idea of the  paper is that our model is not only from the applications point of view, but also from  a new  system of differential equations. That is, we develop differential equations with piecewise constant argument of generalized type  to  a new class of systems; impulsive differential equations with piecewise constant delay and apply them to recurrent neural networks \cite{ak0,ak2,ak3,ak4,ak6,ak7}.  Another novelty  of this paper  is that the sequence of moments $\theta_k,\ k\in \mathbb N,$ where  the  constancy  of the argument changes,  and the sequence of impulsive moments,  $\tau_k,$ are  different.  More precisely,   each moment $\tau_i,\ i\in \mathbb N,$ is an  interior point   of an interval $(\theta_k, \theta_{k+1}).$  This gives  to  our investigations more biological  sense, as well as provides new theoretical  opportunities. 

\section{Model description and preliminaries}\label{model}
Let $\mathbb N=\{0,1,2,\ldots\}$ and $\mathbb R^+=[0,\infty)$  be the sets of  natural and nonnegative real numbers, respectively, and denote a norm on $\mathbb R^m$  by $||\cdot||$ where  $||u||=\displaystyle\sum^{m}_{i=1}|u_{i}|.$ Fix two real valued sequences $\theta=\left\{\theta_k\right\},\tau=\left\{\tau_k\right\},\ k\in \mathbb N,\  \tau\cap\theta=\phi$ such that $\theta_k<\theta_{k+1}$ with $\theta_k\rightarrow \infty$ as $k\rightarrow \infty$  and $\tau_k<\tau_{k+1}$ with $\tau_k\rightarrow \infty$ as $k\rightarrow \infty,$ and there exist two positive numbers $\overline{\theta},\underline{\tau}$ such that $\theta_{k+1}-\theta_{k} \leq \overline{\theta}$ and $\underline{\tau} \leq \tau_{k+1}-\tau_k, \ k\in{\mathbb N}.$  The condition of the empty  intersection is caused by  the investigation reasons. Otherwise,  the proof of auxiliary  results  needs several  additional assumptions.

\begin{figure}[!htp] 
\centering
\includegraphics[width=12cm,height=9cm]{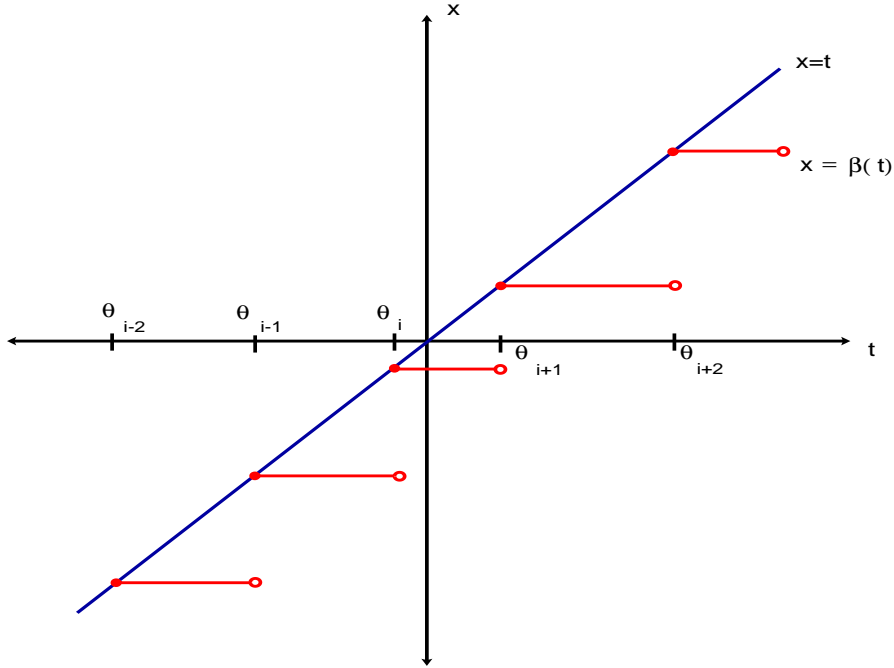}
\caption{The graph of the argument $\beta(t).$}
\label{beta}
\end{figure}

The main subject  under investigation in this paper is the following  impulsive recurrent neural networks  with  piecewise constant delay 
\begin{subequations}\label{ene} 
\begin{eqnarray}
x'_{i}(t) &=& -a_{i}x_{i}(t) + \displaystyle\sum_{j=1}^{m}b_{ij}f_j(x_{j}(t))  \\ &&+\nonumber \sum_{j=1}^{m}c_{ij}g_j(x_{j}(\beta(t)))+d_{i},  \ t \not = \tau_{k}\label{1} \\
\Delta x_{i} \mid_{t=\tau_k}&=& I_{ik}(x_i(\tau_k^-)),\quad a_{i}>0, \ i=1,2,\cdots,m,\quad  k\in \mathbb N, \label{second}
\end{eqnarray}
\end{subequations}
where $ \beta(t)=\theta_{k}$ (see  Fig. \ref{beta}) if  $\theta_{k} \leq t < \theta_{k+1} ,\, k \in{\mathbb N},\,t \in{\mathbb R^+},$  is an identification function, 
$\Delta x_{i}(\tau_k)$ denotes $x_i(\tau_k)-x_i(\tau_k^{-}),$ where $x_i(\tau_k^{-})=\lim_{h\rightarrow 0^{-}}x_i(\tau_k+h).$ Moreover, 
$m$ corresponds to the number of units in a neural
network, $x_{i}(t)$ stands for the state vector of the $i$th unit
at time $t,$ $f_j(x_{j}(t))$ and 
$g_j(x_{j}(\beta(t)))$ denote, respectively, the measures of activation to its incoming 
potentials of the unit $j$ at time $t$  and $\beta(t),$ $b_{ij}, c_{ij}, d_{i}$ are real constants, $b_{ij}$ means the strength of the $j$th
unit on the $i$th unit at time $t,$ $c_{ij}$ infers the strength of the $j$th
unit on the $i$th unit at time $\beta(t),$ $d_{i}$ signifies the external bias
on the $i$th unit and $a_i$ represents the rate with which the $i$th unit will reset its
 potential to the resting state in isolation when it is disconnected
from the network and external inputs.

In the theory of differential equations with piecewise constant argument \cite{ak0,ak2,ak3}, we take the function $\beta(t)=\theta_{k}$ if 
$\theta_{k}\leq t < \theta_{k+1}$, that is, $\beta(t)$ is right continuous. However, as it is usually done in the theory 
of impulsive differential equations, at the points of discontinuity  $\tau_{k}$ of 
the solution, solutions are left continuous. Thus,  the right continuity is more convenient assumption if one considers equations with piecewise constant arguments, and we shall assume  the continuity  for both, impulsive moments and moments of the switching of constancy  of the argument.

We say that the function $\varphi:\mathbb R^+  \rightarrow \mathbb  R^m$ is from the set $PC_{\tau}(\mathbb R^+,{\mathbb R^{m}})$ if:
\begin{enumerate}
 \item [(i)] $\varphi$ is right continuous on $\mathbb R^{+};$
  \item [(ii)] it is continuous everywhere except possibly moments $\tau$ where it has discontinuities of the first kind.
\end{enumerate}

Moreover, we introduce a set of functions $PC_{\tau\cup\theta}(\mathbb R^+,{\mathbb R^{m}})$  if we replace $\tau$ by $\tau\cup\theta$ in the last definition. In our paper, we understand that $\varphi:\mathbb R^+  \rightarrow \mathbb  R^m$ is a solution of (\ref{ene}) if $\varphi \in PC_{\tau}(\mathbb R^+,{\mathbb R^{m}})$ and $\varphi' \in PC_{\tau\cup\theta}(\mathbb R^+,{\mathbb R^{m}}).$

Throughout this paper, we assume the following hypotheses:

\begin{enumerate}
   \item [(H1)] there exist  Lipschitz constants $L_{j}^f,L_{j}^g>0$ such that $$|f_{j}(u)-f_{j}(v)|\leq L_{j}^f|u-v|,$$ 
$$|g_{j}(u)-g_{j}(v)|\leq L_{j}^g|u-v|$$
for all $u,v\in \mathbb R^m,\ j = 1,2,\cdots,m;$
  \item [(H2)] the impulsive operator $I_{ik}: \mathbb R^+\rightarrow \mathbb R^+$ satisfies  $$|I_{ik}(u)-I_{ik}(v)|\leq \ell|u-v|$$
for all $u,v\in \mathbb R^m,\ i = 1,2,\cdots,m,\ k\in \mathbb N$ where $\ell$ is a positive Lipschitz constant.
\end{enumerate}

For the sake of convenience, we adopt the following notations in the sequel:
\begin{eqnarray*}
k_1&=&\displaystyle \max_{1\leq i\leq m}\Big(a_{i}+L_i^f\sum_{j=1}^{m}|b_{ji}|\Big), \ k_2=\displaystyle \max_{1\leq i\leq m}\Big(L_i^g\sum_{j=1}^{m}|c_{ji}|\Big),\ k_3=\displaystyle \max_{k\geq 1}\Big(|I_{k}(0)|\Big),\\ k_4&=&\displaystyle \max_{1\leq i\leq m}\Big(\sum_{j=1}^{m}\Big(|b_{ji}||f_i(0)|+ |c_{ji}||g_i(0)|\Big)\Big).
\end{eqnarray*}

Denote by $p_{k}$ the number of points $\tau_{i}$ in the interval $(\theta_{k}, \theta_{k+1}),\ k\in \mathbb N.$  We assume that   $p=\displaystyle \max_{k\in \mathbb N} p_{k} < \infty.$

Assume, additionally, that 
\begin{enumerate}
 \item [(H3)] $\Big[(k_1+2k_2)\overline{\theta}+\ell p\Big]\Big(1+\ell\Big)^{p}\e^{k_1\overline{\theta}}<1;$
 \item [(H4)] $k_2\overline{\theta}+(k_1\overline{\theta}+\ell p)(1+k_2\overline{\theta})(1+\ell)^{p}\e^{k_1\overline{\theta}}<1.$
\end{enumerate}

\subsection{Equilibrium for neural networks models with discontinuities}

In our paper, the main objects of discussion are the concepts of equilibria for discontinuous neural networks models. 

First of all, we must say that an equilibrium is a necessary attribute of any models for neural networks. Since, in the applications of Hopfield-type neural networks, when these networks can be considered as a nonlinear associative memory, the equilibrium states of the network serve as stored patterns and their stability means that the stored patterns can be retrieved representing the fundamental memories of the system \cite{a115,a1192,hop,a121,HNNWI1,HNNWI3,HNNWI4,HNNWDI1,HNNWDI2,HNNWDI3}. To give more sense, the stable points of the phase space of the network are the fundamental memories, or prototype states of the network. For example, when a network has a pattern containing partial but sufficient information about one of the fundamental memories, we may represent it as a starting point in the phase space. Provided that the starting point is close to the stable point representing the memory being retrieved, finally, the system converges onto the memory state itself. Consequently, it can be said that Hopfield network is a dynamical system whose phase space contains a set of fixed (stable) points representing the fundamental memories of the system \cite{a115,a1192,hop,a121}. This is obviously  true for  networks, which  generalize the Hopfield one. For example, Recurrent neural networks, which  are  considered in our manuscript, Cohen-Grossberg neural Networks (Hopfield neural networks as  a special version) \cite{a114} and Cellular neural networks \cite{CNNT}.

Now,  we want  to  give a general  analysis of  the equilibrium  concept  for discontinuous dynamics descibed by  impulsive differential  equations. Consider a  $x(t),$ motion of a neural network,  assuming  that  it  is discontinuous. Since the process is assumed right continuous, we have that $x(\tau_k)-x(\tau_k^{-})= I_k(x(\tau_k^-)),\ k\in \mathbb N,$ where $I_k(x)$ is an activation function. If one expects that $x=x^*$ is an equilibria of the network, the necessarily $I_{k}(x^*)=0,$ since $x=x^*$ is not a fixed point, otherwise.  Thus,  this condition is to be the cornerstone of our manuscript.  One can  see  that  the property  of impulses to  diminish  at  the equilibrium is not  easy  to   derive  from the mathematical  analysis as it is true for  continuous dynamics \cite{mfp,cgp}, and  it  should be  introduced in a special  way \cite{HNNWI1,HNNWI3,HNNWI4,HNNWDI1,HNNWDI2,HNNWDI3}. In  our paper we  pay attention to the main subject more than usually. One must  emphasise that despite zero  impulses at  equilibria, motions near eqiulibria admit  non-zero  impulses, and it  makes the systems considerable in the dynamics theory. Thus, if  one  recognize  the zero  right-hand side of  a differential equation  at  an  equilibrium, then   annihilated  impulses  have to  be accepted, too. Mechanically, the zero  impulses are the same as the zero  velocity. If motionless points recognized, they  must  be  accepted not  only  for smooth systems but also for  discontinuous systems. To  illustrate  this point of view let us remember  the impulsive differential equation  in its general  form 
\begin{eqnarray*}
\begin{array}{l}
 x' = f(t,x),  \ t \not = t_{k}\label{1aa} \\
\Delta x \mid_{t=t_k}= I_k(x).
\end{array}
\end{eqnarray*}
The system can  be written as 
\begin{eqnarray*}
x' =f(t,x)+ \sum^{\infty}_{k=1}I_k(x)\delta(t-t_k), \label{2aa} 
\end{eqnarray*}
where $\delta(t)$ is the   delta function. Then the equation for  equilibria has  the form  
\begin{eqnarray}\label{3aa}
f(t,x^*)+ \sum^{\infty}_{k=1}I_k(x^*)\delta(t-t_k)=0, 
\end{eqnarray}
which  is very  similar to  that  for  ordinary  differential  equations.

It  is obvious that   equilibria are  easier  for analysis  if  equations are  autonomous, but    they  are also  very  popular for non-autonomous systems,  in theory  as well  as in applications.  Investigation  of   eqution (\ref{3aa}) is not  more difficult, than that   of  $f(t,x) =0,$ in general.  We suppose that   theoretically  proven   existence of equilibria give  us possibility  to  investigate concrete examples numerically and by  simulations,  admitting  small errors,  as it  is done in our  example with  equation  (\ref{ex3}) in Section  \ref{simulations}.

We denote  the constant vector by $x^*=(x_{1}^*,\ldots,x_{m}^*)^T\in \mathbb R^m$, where  the components $x_{i}^*$ are governed by the algebraic system
\begin{eqnarray}\label{eqq}
0=-a_{i}x_{i}^*+\displaystyle\sum_{j=1}^{m}b_{ij}f_j(x_{j}^*) + \sum_{j=1}^{m}c_{ij}g_j(x_{j}^*))+d_{i}.
\end{eqnarray}

%We denote an equilibrium state for the differential equation of (\ref{1}) by the constant vector $x^*=(x_{1}^*,\ldots,x_{m}^*)^T\in \mathbb R^m$, where  the components $x_{i}^*$ are governed by the algebraic system
%\begin{eqnarray*}
%0=-a_{i}x_{i}^*+\displaystyle\sum_{j=1}^{m}b_{ij}f_j(x_{j}^*) + \sum_{j=1}^{m}c_{ij}g_j(x_{j}^*))+d_{i}.
%\end{eqnarray*}

The proof of following lemma  is very similar to that of  Lemma 2.2 in \cite{HNNWI1} and therefore we omit it here. 

\begin{lemma}\label{eueq}
Assume $\rm{(H1)}$ holds. If  the condition  
\begin{eqnarray*}
a_{i}>L_{i}^f\displaystyle\sum_{j=1}^{m}|b_{ji}|+L_{i}^g\displaystyle\sum_{j=1}^{m}|c_{ji}| ,\quad i=1,\ldots,m.
\end{eqnarray*}
is satisfied, then there exists  a  unique  constant vector $x^*=(x_{1}^*,\ldots,x_{m}^*)^T$ such that (\ref{eqq}) is valid.
\end{lemma}
In this case, $x^*=(x_{1}^*,\ldots,x_{m}^*)^T$ is an equilibrium point of equation (\ref{1}).

Let us denote the set of all zeros  of the impact activation functions by $$\Omega=\left\{x=(x_{1},\ldots,x_{m})^T\in \mathbb R^m \mid I_{ik}(x)=0,\ i=1,\ldots,m,\ k\in \mathbb N \right\}.$$
From now on, we shall need the following assumption:
\begin{enumerate}
 \item [(A)] $x^{*} \in \Omega.$
\end{enumerate}

%\begin{theorem}
%If the equilibrium $x^*=(x_{1}^*,\ldots,x_{m}^*)^T\in \mathbb R^m$ of the differential equation of (\ref{1}) satisfies $I_{k}(x_{i}^*)=0$ for all $i=1,\ldots,m,\ k\in \mathbb N.$ Then,  $x^*$ is an equilibrium point of (\ref{1}). 
%\end{theorem}

\begin{theorem}
If $x^*=(x_{1}^*,\ldots,x_{m}^*)^T$  is an equilibrium point of  equation (\ref{1}) and the condition $\rm{(A)}$ is satisfied, then  $x^*$ is an equilibrium point of (\ref{ene}).
\end{theorem}

Now we need the following equivalence lemmas which will be used in the proof of next assertions. The proofs are omitted here, since it is similar  in \cite{a0000,ak0,ak2,ak3,ak6,SP}.

\begin{lemma}\label{irep0}
A function $x(t)=x(t,t_{0},x^{0}) = (x_1(t),\cdots,x_m(t))^T,$ where $t_{0}$ is a fixed real  number, is a solution of (\ref{ene}) on $\mathbb R^{+}$ 
if and only if it is a solution, on $\mathbb R^{+},$ of the following integral equation: 
\begin{eqnarray*}
x_i(t) &=& e^{-a_i(t-t_{0})}x_i^{0} + \displaystyle\int^{t}_{t_{0}}e^{-a_i(t-s)} \left[\sum_{j=1}^{m}b_{ij}f_j(x_j(s))\right. \nonumber \\&&
+\left. \displaystyle\sum_{j=1}^{m}c_{ij}g_j(x_j(\beta(s)))+d_i\right]ds + \displaystyle\sum_{t_{0}\leq \tau_{k}<t}e^{-a_i(t-\tau_{k})}I_{ik}(x_i(\tau_{k}^{-})),
\end{eqnarray*}
for $i=1,\cdots,m,\ t\geq t_0.$
\end{lemma}

\begin{lemma}\label{irep}
A function $x(t)=x(t,t_{0},x^{0})= (x_1(t),\cdots,x_m(t))^T,$ where $t_{0}$ is a fixed real  number, is a solution of (\ref{ene}) on $\mathbb R^{+}$ 
if and only if it is a solution, on $\mathbb R^{+},$ of the following integral equation: 
\begin{eqnarray*}
x_i(t) &=& x_i^{0} + \displaystyle\int^{t}_{t_{0}} \left[-a_ix_i(s)+\sum_{j=1}^{m}b_{ij}f_j(x_j(s))\right. \nonumber \\&&
+\left. \displaystyle\sum_{j=1}^{m}c_{ij}g_j(x_j(\beta(s)))+d_i\right]ds + \displaystyle\sum_{t_{0}\leq \tau_{k}<t}I_{ik}(x_i(\tau_{k}^{-})),
\end{eqnarray*}
for $i=1,\cdots,m,\ t\geq t_0.$
\end{lemma}

Consider the following system

\begin{eqnarray}
\left\{\begin{array}{ll}
 x'_{i}(t) = -a_{i}x_{i}(t) + \displaystyle\sum_{j=1}^{m}b_{ij}f_j(x_{j}(t)) + \sum_{j=1}^{m}c_{ij}g_j(x_{j}(\theta_r))+d_{i},  \ t \not = \tau_{r}\label{s1} \\
\Delta x_{i} \mid_{t=\tau_r}= I_{ir}(x_i(\tau_r^-)),\quad  i=1,2,\cdots,m,\quad  k\in \mathbb N.
\end{array}\right.
\end{eqnarray}

In the next lemma the conditions of existence and uniqueness of solutions  are established for arbitrary initial moment $\xi.$

\begin{lemma}\label{eu}
Assume that conditions $\rm{(H1)-(H3)}$ are fulfilled, and fix $r \in \mathbb N.$ Then for every $(\xi,x^{0})\in [\theta_r,\theta_{r+1}]\times\mathbb R^m$  there exists a unique solution 
$x(t)=x(t, \xi, x^{0})=(x_1(t),\cdots,x_m(t))^T$
of  (\ref{s1}) on $[\theta_r,\theta_{r+1}]$  with $x(\xi)=x^{0}.$
\end{lemma}

\noindent {\bf Proof.} $Existence:$ 
Denote $ \vartheta(t)=x(t,\xi,x^{0}),\vartheta(t)=(\vartheta_1(t),\cdots, \vartheta_m(t))^T.$ From Lemma \ref{irep}, we have
\begin{eqnarray}
\vartheta_i(t) &=& x_i^{0} + \displaystyle\int^{t}_{\xi}\left[-a_{i}\vartheta_{i}(s)+ \sum_{j=1}^{m}b_{ij}f_j(\vartheta_j(s)) 
+ \displaystyle\sum_{j=1}^{m}c_{ij}g_j(\vartheta_j(\theta_r))+d_i\right]ds \nonumber \\&& + \displaystyle\sum_{\xi\leq \tau_{r}<t}I_{ir}(\vartheta_i(\tau_{r}^{-})).\label{6}
\end{eqnarray}
Define a norm  $||\vartheta(t)||_0=\displaystyle\max_{[\theta_{r},\ \theta_{r+1}]}||\vartheta(t)||$ and construct the following sequences  $\vartheta_i^{n}(t),\ \vartheta_{i}^0(t) \equiv x_i^{0},\ i = 1,\ldots,m, \ n\geq 0$ such that

\begin{eqnarray*}
\vartheta_i^{n+1}(t) &=& x_i^{0} + \displaystyle\int^{t}_{\xi}\left[-a_{i}\vartheta_{i}^{n}(s)+ \sum_{j=1}^{m}b_{ij}f_j(\vartheta_j^{n}(s)) 
+ \displaystyle\sum_{j=1}^{m}c_{ij}g_j(\vartheta_j^{n}(\theta_r))+d_i\right]ds \nonumber \\&& + \displaystyle\sum_{\xi\leq \tau_{r}<t}I_{ir}(\vartheta_i^{n}(\tau_{r}^{-})).
\end{eqnarray*}

One can find that
$$||\vartheta^{n+1}(t)-\vartheta^{n}(t)||_0\leq\left(\left[(k_1+k_2)\overline{\theta}+\ell p\right]\right)^{n}\kappa,$$ where $$\kappa=\left(\left[(k_1+k_2)\overline{\theta}+\ell p\right]||\vartheta^{0}||+\overline{\theta}\Big(\displaystyle\sum_{i=1}^{m}d_i\Big)+\overline{\theta}mk_4+mpk_{3}\right).$$
Since the condition $\rm{(H3)}$ implies $\left[(k_1+k_2)\overline{\theta}+\ell p\right]< 1,$ then  the sequences $\vartheta_i^{n}(t)$ are convergent and their limits satisfy (\ref{6}) on 
$[\theta_{r},\ \theta_{r+1}].$   The existence is proved.

$Uniqueness:$ Let us denote the  solutions of (\ref{ene}) by  $x^1(t)=x(t,\xi,x^1),$ $x^2(t)=x(t,\xi,x^2),$  where $\theta_r \leq \xi \leq \theta_{r+1}.$
It is sufficient to check that for each $t\in [\theta_{r},\theta_{r+1}],$ and  $x^2=(x_{1}^2,\cdots,x_{m}^2)^T, x^1=(x_{1}^1,\cdots,x_{m}^1)^T \in \mathbb R^m, x^2\neq x^1,$ the condition $x^1(t)\neq x^2(t).$ Then, we have

\begin{eqnarray*}
||x^1(t)-x^2(t)||&\leq& ||x^1-x^2||+\sum_{i=1}^{m}\left\{\int_{\xi}^{t}\left[\Big(a_{i}+
\sum_{j=1}^{m}L_{i}^{f}|b_{ji}|\Big)\left|x_{i}^2(s)-x_{i}^1(s)\right|\right.\right.\nonumber\\&&+\left.\left.\sum_{j=1}^{m}L_{i}^{g}|c_{ji}| \left|x_{i}^2(\theta_r)-x_{i}^1(\theta_r)\right|\right]ds \nonumber + \ell \displaystyle\sum_{\xi\leq\tau_{r}<t}\left|x_i^{2}(\tau_{r}^{-})-x_i^{1}(\tau_{r}^{-})\right|
\right\}\\&\leq&||x^1-x^2||+k_2\overline{\theta}||x^1(\theta_r)-x^2(\theta_r)||+ k_1\int_{\xi}^{t}||x^1(s)-x^2(s)||ds\\&&+\ell \displaystyle\sum_{\xi\leq\tau_{r}<t}||x^1(\tau_{r}^{-})-x^2(\tau_{r}^{-})||. 
\end{eqnarray*}

Using Gronwall-Bellman Lemma for piecewise continuous functions \cite{SP,laks},  one can  obtain that
\begin{eqnarray*}
||x^1(t)-x^2(t)||\leq \Big( ||x^1-x^2||+k_2\overline{\theta}||x^1(\theta_r)-x^2(\theta_r)||\Big)(1+\ell)^{p}\e^{k_1\overline{\theta}}.
\end{eqnarray*}
Particularly, 
\begin{eqnarray*}
||x^1(\theta_r)-x^2(\theta_r)||\leq \Big( ||x^1-x^2||+k_2\overline{\theta}||x^1(\theta_r)-x^2(\theta_r)||\Big)(1+\ell)^{p}\e^{k_1\overline{\theta}}.
\end{eqnarray*}
Hence,
\begin{eqnarray}\label{inq1}
||x^1(t)-x^2(t)||\leq \left[\frac{(1+\ell)^{p}\e^{k_1\overline{\theta}}}{1-k_2\overline{\theta}(1+\ell)^{p}\e^{k_1\overline{\theta}}}\right] ||x^1-x^2||.
\end{eqnarray}
Also, we peculiarly have

\begin{eqnarray}\label{inq2}
||x^1(\tau_{r}^{-})-x^2(\tau_{r}^{-})||\leq \left[\frac{(1+\ell)^{p}\e^{k_1\overline{\theta}}}{1-k_2\overline{\theta}(1+\ell)^{p}\e^{k_1\overline{\theta}}}\right] ||x^1-x^2||.
\end{eqnarray}
 
On the other hand, assume on the contrary that there exists $t\in [\theta_r,\theta_{r+1}]$
such that $x^1(t)=x^2(t).$ Then

\begin{eqnarray}
||x^1-x^2||&\leq& \sum_{i=1}^{m}\left\{\int_{\xi}^{t}\left[\Big(a_{i}+
\sum_{j=1}^{m}L_{i}^{f}|b_{ji}|\Big)\left|x_{i}^2(s)-x_{i}^1(s)\right|\right.\right.\nonumber\\&&+\left.\left.\sum_{j=1}^{m}L_{i}^{g}|c_{ji}| \left|x_{i}^2(\theta_r)-x_{i}^1(\theta_r)\right|\right]ds \nonumber + \ell \displaystyle\sum_{\xi\leq\tau_{r}<t}\left|x_i^{2}(\tau_{r}^{-})-x_i^{1}(\tau_{r}^{-})\right|
\right\}\\&\leq&k_1\int_{\xi}^{t}||x^1(s)-x^2(s)||ds+k_2\overline{\theta}||x^1(\theta_r)-x^2(\theta_r)|| \nonumber \\&&+\ell p||x^1(\tau_{r}^{-})-x^2(\tau_{r}^{-})||. \label{inq3}
\end{eqnarray}

Consequently, substituting (\ref{inq1}) and (\ref{inq2}) in (\ref{inq3}), we obtain

\begin{eqnarray}\label{inq4}
||x^1-x^2||\leq\Big[(k_1+2k_2)\overline{\theta}+\ell p\Big]\Big(1+\ell\Big)^{p}\e^{k_1\overline{\theta}}||x^1-x^2||.
\end{eqnarray}

Thus, one can see that (H3) contradicts with (\ref{inq4}). The lemma is proved. $\Box$

\begin{theorem}\label{eus}
Assume that conditions $\rm{(H1)-(H3)}$ are fulfilled. Then, for every $(t_{0},x^{0})\in \mathbb R^+ \times \mathbb R^m,$ there
exists a unique solution $x(t)=x(t,t_{0},x^{0})=(x_1(t),\cdots,x_m(t))^T,\ t \geq t_{0},$ of (\ref{ene}), such
that $x(t_{0})=x^{0}.$
\end{theorem}

\noindent {\bf Proof.} 
Fix $t_{0}\in \mathbb R^+.$ It is clear that there exists $r\in \mathbb N$ such that $t_{0}\in[\theta_{r},\theta_{r+1}).$ 
Using previous lemma for $\xi=t_{0},$ one can obtain  that there exists a unique solution
$x(t)=x(t,t_{0},x^{0})$  on $[\xi,\theta_{r+1}].$ Next, we again apply the last lemma to obtain the unique  solution on interval  $[\theta_{r+1},\theta_{r+2}).$ The mathematical induction completes the proof. $\Box$

\section{Global asymptotic stability }\label{gas}
In this section, we will focus our attention on  giving  sufficient conditions for the global asymptotic stability of the equilibrium, $x^*,$ of (\ref{ene}) based on linearization \cite{ak3,SP}.

The system (\ref{ene}) can be simplified as follows. Substituting $y(t)=x(t)-x^*$ into (\ref{ene}) leads to 
\begin{eqnarray}
\left\{\begin{array}{ll}
y'_{i}(t) = -a_{i}y_{i}(t) + \displaystyle\sum_{j=1}^{m}b_{ij}\phi_j(y_{j}(t)) + \sum_{j=1}^{m}c_{ij}\psi_j(y_{j}(\beta(t))), \ t\neq \tau_{k} \label{11} \\
\Delta y_{i} \mid_{t=\tau_k}= W_{ik}(y_i(\tau_k^-)),\quad i=1,2,\cdots,m,\quad  k\in \mathbb N,
\end{array}\right.
\end{eqnarray}
where 
$\phi_{j}(y_{j}(t))=f_{j}(y_{j}(t)+x_{j}^*)-f_{j}(x_{j}^*),\   
\psi_{j}(y_{j}(t))=g_{j}(y_{j}(t)+x_{j}^*)-g_{j}(x_{i}^*)$ and 
$W_{ik}(y_i(\tau_k^-))=I_{ik}(y_i(\tau_k^-)+x_{i}^*)-I_{ik}(x_{i}^*).$  From hypotheses (H1) and (H2), we have the following inequalities: $|\phi_{j}(\cdot)|\leq L_j^{f}|(\cdot)|,\ |\psi_{j}(\cdot)|\leq L_j^{g}|(\cdot)|$ and $|W_{ik}(\cdot)|\leq \ell |(\cdot)|.$ 

It is clear that the stability of the zero solution of  (\ref{11}) is equivalent to the
stability of the equilibrium $x^*$ of  (\ref{ene}). Therefore, in what follows, we  discuss the stability of the zero solution of (\ref{11}).

First of all, we give the  lemma below which is one of the most  important results of the present paper. One can see that this lemma is generalized version of the lemmas   in \cite{ak0,ak2,ak3,ak4,ak6,ak7}.

For simplicity of notation, we denote
$$\lambda=\left(1-\Big(k_2\overline{\theta}+(k_1\overline{\theta}+\ell p)(1+k_2\overline{\theta})(1+\ell)^{p}\e^{k_1\overline{\theta}}\Big)\right)^{-1}.$$

\begin{lemma} 
Let $y(t)=(y_1(t),\cdots,y_m(t))^T$ be a solution of (\ref{11}) and $\rm{(H1)-(H4)}$ be satisfied. Then, the following inequality 
\begin{eqnarray}
||y(\beta(t))||\leq  \lambda||y(t)|| \label{2}
\end{eqnarray}
holds for all $t\in \mathbb R^+.$ 
\end{lemma}

\noindent {\bf Proof.} 
Fix $t\in \mathbb R^+,$ there exists  $k\in \mathbb N$ such that 
$t\in [\theta_{k},\theta_{k+1}).$ Then, from Lemma \ref{irep}, we have
\begin{eqnarray*}
||y(t)||&=&\sum_{i=1}^{m}|y_i(t)|\\&\leq& ||y(\theta_{k})||+\sum_{i=1}^{m}\left\{\int_{\theta_{k}}^{t}\left[\Big(a_{i}+
\sum_{j=1}^{m}L_{i}^{f}|b_{ji}|\Big)\left|y_{i}(s)\right|\right.\right.\nonumber\\&&+\left.\left.\sum_{j=1}^{m}L_{i}^{g}|c_{ji}| \left|y_{i}(\theta_k)\right|\right]ds\right.+\left.\ell\displaystyle\sum_{t_0\leq\tau_{k}<t}\left|y_i(\tau_{k}^{-})\right|
\right\}\\&\leq&(1+k_2\overline{\theta})||y(\theta_{k})||+k_1\int_{\theta_{k}}^{t}||y(s)||ds+\ell \displaystyle\sum_{t_0\leq\tau_{k}<t}||y(\tau_{k}^{-})||.
\end{eqnarray*}
Applying the analogue of  Gronwall-Bellman Lemma \cite{SP,laks}, we obtain
\begin{eqnarray}
||y(t)||\leq (1+k_2\overline{\theta})(1+\ell)^{p}\e^{k_1\overline{\theta}}||y(\theta_{k})||. \label{4}
\end{eqnarray}
Particularly,
\begin{eqnarray}
||y(\tau_k^{-})||\leq(1+k_2\overline{\theta})(1+\ell)^{p}\e^{k_1\overline{\theta}} ||y(\theta_{k})||.\label{tau1}
\end{eqnarray}
Moreover, for $t\in [\theta_{k},\theta_{k+1}),$ we also have
\begin{eqnarray*}
||y(\theta_{k})||&\leq&||y(t)||+k_2\overline{\theta}||y(\theta_{k})||+k_1\int_{\theta_{k}}^{t}||y(s)||ds\\&&+\ell\displaystyle\sum_{t_0\leq \tau_{k}<t}||y(\tau_{k}^{-})||. 
\end{eqnarray*}
The last inequality together with   (\ref{4}) and (\ref{tau1}) imply
\begin{eqnarray*}
||y(\theta_{k}))||&\leq& ||y(t)||+ \Big[k_2\overline{\theta}+(k_1\overline{\theta}+\ell p)(1+k_2\overline{\theta})(1+\ell)^{p}\e^{k_1\overline{\theta}}\Big]||y(\theta_{k})||
\end{eqnarray*}
Thus, we have from condition  (H4) that
$$||y(\theta_{k})||\leq \lambda||y(t)||,\,\quad t\in[\theta_{k},\theta_{k+1}).$$
Therefore, (\ref{2}) holds  for all $t\in \mathbb R^+.$ This completes the  proof of lemma. $\Box$

Now, we are ready to give sufficient conditions for the global asymptotic 
stability of (\ref{ene}). For convenience, we adopt the notation given below in the sequel:
$$\mu=\displaystyle \max_{1\leq i\leq m}\Big(L_i^f\sum_{j=1}^{n}|b_{ji}|\Big).$$

From now on we need  the following  assumption:
\begin{enumerate}
\item[(H5)]$\gamma-\mu-\lambda k_2-\frac{ln(1+l)}{\underline{\tau}}>0,\quad \gamma=\displaystyle \min_{1\leq i\leq m}a_{i}.$
\end{enumerate}

The next theorem is a modified version of the theorem in \cite{SP}, for our system.
\begin{theorem}\label{tgs}
Assume  that $\rm{(H1)-(H5)}$ are fulfilled.Then, 
the zero solution of (\ref{11}) is globally asymptotically stable.
\end{theorem}

\noindent {\bf Proof.}
Let $y(t)=(y_1(t),\cdots,y_m(t))^T$ be an arbitrary solution of (\ref{11}). From Lemma \ref{irep0}, we have
\begin{eqnarray*}
||y(t)|| &\leq&e^{-\gamma(t-t_0)}||y_0|| +\sum_{i=1}^{m}\left\{ \displaystyle\int^{t}_{t_{0}}e^{-\gamma(t-s)}\left[\displaystyle\sum_{j=1}^{m}L_i^{f}|b_{ji}||y_i(s)| \right.\right.\\&&+ \left.\left.\sum_{j=1}^{m}L_i^{g}|c_{ji}||y_i(\beta(s))|\right]ds
+\ell\displaystyle\sum_{t_0\leq\tau_{k}<t}e^{-\gamma(t-\tau_k)}|y_i(\tau_{k}^{-})|\right\}\\&\leq&e^{-\gamma(t-t_0)}||y_0||+\left(\mu+\lambda k_2 \right)\displaystyle\int^{t}_{t_{0}}e^{-\gamma(t-s)}||y(s)||ds\\&&+\ell \displaystyle\sum_{t_0\leq\tau_{k}<t}e^{-\gamma(t-\tau_k)}||y(\tau_{k}^{-})||.
\end{eqnarray*}
Then, we can write the last  inequality as,
\begin{eqnarray*}
e^{\gamma(t-t_{0})}||y(t)|| &\leq& ||y_0||+\left(\mu+\lambda k_2 \right)\displaystyle\int^{t}_{t_{0}}e^{\gamma(s-t_0)}||y(s)||ds\\&&+\ell\displaystyle\sum_{t_0\leq\tau_{k}<t}e^{\gamma(\tau_k-t_0)}||y(\tau_{k}^{-})||.
\end{eqnarray*}
By virtue of  Gronwall-Bellman Lemma \cite{SP}, we obtain
\begin{eqnarray*}
e^{\gamma(t-t_{0})}||y(t)||\leq e^{\left(\mu+\lambda k_2 \right)(t-t_{0})}[1+\ell]^{i(t_{0},t)}||y_0||,
\end{eqnarray*}
where $i(t_0, t)$ is the number of points $\tau_k$ in $[t_{0}, t).$ Then, we have
$$||y(t)||\leq  e^{-(\gamma-\mu-\lambda k_2-\frac{ln(1+\ell)}{\underline{\tau}})(t-t_{0})}||y_{0}||. $$
Hence, using  the condition (H5), we see that the zero solution of system (\ref{11}) is globally asymptotically stable. $\Box$

\section{Existence of periodic solutions}\label{seceql}

In this section,  we shall discuss the  existence of periodic solution of  (\ref{ene}) and its stability. To do so, we  need the following assumptions:
\begin{enumerate}
\item[(H6)] the sequences $\tau_{k}$ and  $\theta_k,\ k\in \mathbb N$ satisfy  $(\omega, p)$ and  $(\omega, p_{1})$-properties; that is, there are positive integers $p$ and  $p_1$ such that the equations $\tau_{k+p}=\tau_{k}+\omega$ and $\theta_{k+p_{1}}=\theta_{k}+\omega$ hold for all $k\in \mathbb N$
and $I_{k+p}=I_{k}$ for a fixed positive real period $\omega.$ 
\item[(H7)] $\alpha_1=\mathcal R\Big(\omega\left(\mu+\lambda k_2\right)+\ell p\Big)<1,$\quad where $\mathcal R=\frac{1}{1-\e^{-\gamma\omega}}.$ 
\end{enumerate}
For $\tau_k$ and $\theta_{k},$  let $[0,\omega]\cap\{\tau_k\}_{k\in \mathbb N}=\{\tau_1,\cdots,\tau_p\}$ and $[0,\omega]\cap\{\theta_k\}_{k\in \mathbb N}=\{\theta_1,\cdots,\theta_{p_{1}}\},$ respectively.

Here, we will give the following version of the Poincare' criterion for system (\ref{ene}) which  can be easily proved (see, also, \cite{SP}).

\begin{lemma}
Suppose that conditions $\rm{(H1)-(H3)}$ and $\rm{(H6)}$ are valid. Then, solution $x(t)=x(t,t_{0},x^{0})=(x_1,\cdots,x_m)^T$ of (\ref{ene}) with $x(t_{0})=x^{0}$
is \\$\omega-$periodic if and only if \ $x(\omega)=x(0).$
\end{lemma}

\begin{theorem}\label{ups}
Assume that conditions $\rm{(H1)-(H3)}$ and  $\rm{(H6)- (H7)}$  are valid. 
 Then system (\ref{ene}) has a unique $\omega-$periodic solution.
\end{theorem}

\noindent {\bf Proof.}
To begin with, let us  introduce a Banach space of periodic functions   $PC_{\omega}=\{\varphi\in PC_{\tau\cup\theta}(\mathbb R^{+}, \mathbb R^m) \mid  \varphi(t+\omega)=\varphi(t),\ t\geq 0 \}$  with the norm $||\varphi||_0=\displaystyle\max_{0\leq t\leq \omega}||\varphi(t)||.$ 

Let $\varphi(t)=(\varphi_{1}(t),\cdots,\varphi_{m}(t))^T \in PC_{\omega}$  satisfying  the inequality $||\varphi(t)||_{0}\leq h.$   Using Lemma 2.1, similarly to the proof  in \cite{SP}, one can show that if $\varphi \in PC_{\omega}$ then the system
\begin{eqnarray*}
\left\{\begin{array}{l}
 x_{i}'(t) = -a_{i}x_{i}(t) + \displaystyle\sum_{j=1}^{m}b_{ij}f_j(\varphi_j(t)) + \sum_{j=1}^{m}c_{ij}g_j(\varphi_j(\beta(t)))+d_{i}, \ t\neq \tau_{k}, \\
\Delta x_{i} \mid_{t=\tau_k}= I_{ik}(\varphi_i(\tau_k^-)),\ i=1,\cdots,m,\quad k=1,2,\cdots,p
\end{array}\right.
\end{eqnarray*}
has the unique $\omega-$ periodic solution
\begin{eqnarray*} 
x_{i}^*(t)&=&\int_0^{\omega}
\mathcal G_{i}(t,s)\left[\displaystyle\sum_{j=1}^{m}b_{ij}f_j(\varphi_j(s)) + \sum_{j=1}^{m}c_{ij}g_j(\varphi_j(\beta(s)))+d_{i}\right]ds \\&&+\sum_{k=1}^{p}
\mathcal G_i(t,\tau_{k})I_{ik}(\varphi_i(\tau_k^{-})), 
\end{eqnarray*}
where 
\begin{eqnarray*}
\mathcal G_{i}(t,s)=(1-e^{-a_{i}\omega})^{-1}\left\{\begin{array}{ll}e^{-a_{i}(t-s)},
& 0\leq s\leq t\leq  \omega, \\
e^{-a_{i}(\omega+t-s)}, & 0\leq
t<s\leq \omega,
\end{array}\right. 
\end{eqnarray*}

which is known as  {\it Green's function} \cite{SP}. Then, one can easily find that\\ $\displaystyle\max_{t,s\in [0,\omega]}\left|\{\mathcal G_{i}(t,s)\}_{i=1,\cdots,m}\right|=\frac{1}{1-\e^{-a_i\omega}}.$

Define the operator $ \mathcal F: PC_{\omega}\rightarrow PC_{\omega}$ such that if $\varphi \in PC_{\omega},$ then
\begin{eqnarray*}
(\mathcal F\varphi)_i(t)&=&\int_0^{\omega}
\mathcal G_{i}(t,s)\left[\displaystyle\sum_{j=1}^{m}b_{ij}f_j(\varphi_j(s)) + \sum_{j=1}^{m}c_{ij}g_j(\varphi_j(\beta(s)))+d_{i}\right]ds \\&&+\sum_{k=1}^{p}
\mathcal G_i(t,\tau_{k})I_{ik}(\varphi_i(\tau_k^{-})), \quad i=1,\cdots,m. 
\end{eqnarray*}
Now, we need to prove that $\mathcal F$ maps $PC_{\omega}$ into itself. That is, we shall show that $\mathcal F \varphi \in PC_{\omega}$ for any $\varphi\in PC_{\omega}.$ It is easy to check that $(\mathcal F\varphi)(t)=((\mathcal F\varphi)_1,\cdots,(\mathcal F\varphi)_m)^T$ is $\omega-$periodic function. Now, if  $\varphi\in PC_{\omega},$ then 
\begin{eqnarray*}
||{\mathcal F}\varphi||&=&\sum^{m}_{i=1}\left|\int_0^{\omega}\mathcal G_{i}(t,s) \left[\displaystyle\sum_{j=1}^{m}b_{ij}f_j(\varphi_j(s))+ \sum_{j=1}^{m}c_{ij}g_j(\varphi_j(\beta(s)))+d_{i}\right]ds \right.\\ 
&&\left.+\sum_{k=1}^{p}\mathcal G_{i}(t,\tau_{k})I_{ik}(\varphi_i(\tau_k^{-}))\right|\\&\leq& \sum^{m}_{i=1}\frac{1}{1-\e^{-a_i\omega}}\left\{\int_0^{\omega}\left[\sum_{j=1}^{m}L_{j}^{f}|b_{ij}||\varphi_j(s)|+ \sum_{j=1}^{m} L_{j}^{g}|c_{ij}||\varphi_j(\beta(s))|\right.\right.\\&&+\left.\left.\sum_{j=1}^{m}|b_{ij}||\varphi_j(0)|+\sum_{j=1}^{m}|c_{ij}||\varphi_j(0)|+d_{i}\right]ds+\sum_{k=1}^{p}|I_{ik}(\varphi_i(\tau_k^{-}))|\right\}\\&\leq&\mathcal R\sum^{m}_{i=1}\left\{\int_0^{\omega}\left[\sum_{j=1}^{m}L_{i}^{f}|b_{ji}||\varphi_i(s)|+ \sum_{j=1}^{m} L_{i}^{g}|c_{ji}||\varphi_i(\beta(s))|\right.\right.\\&&+\left.\left.\sum_{j=1}^{m}|b_{ji}||\varphi_i(0)|+\sum_{j=1}^{m}|c_{ji}||\varphi_i(0)|+d_{i}\right]ds+\ell\sum_{k=1}^{p}|\varphi_i(\tau_k^{-})|+\sum_{k=1}^{p}|I_{ik}(0)|\right\}\\&\leq& \mathcal R\left(\int_0^{\omega}\Big[\mu||\varphi(s)||+k_2||\varphi(\beta(s))||\Big]ds+\ell\sum_{k=1}^{p}||\varphi(\tau_k^{-})||+\omega\Big(\displaystyle\sum^{m}_{i=1}d_{i}\Big)+\omega mk_4+mpk_{3}\right).
\end{eqnarray*}
In this periodical case, we take $k_3=\displaystyle \max_{1\leq k\leq p}\Big(|I_{k}(0)|\Big).$ Thus, it follows that
\begin{eqnarray*}
||{\mathcal F}\varphi||_0&\leq&\mathcal R\Big(\Big(\omega\left(\mu+\lambda k_2\right)+\ell p\Big)||\varphi||_0+ \omega\Big(\displaystyle\sum^{m}_{i=1}d_{i}\Big)+\omega mk_4+mpk_{3}\Big)\\&\leq&\alpha_{1}h+\alpha_{2}.
\end{eqnarray*}
Choose $h$ such that $\alpha_{2}\leq h(1-\alpha_{1}),$
where $\alpha_{2}=\mathcal R \Big(\omega\Big(\displaystyle\sum^{m}_{i=1}d_{i}\Big)+\omega mk_4+mpk_{3}\Big).$ 
Then, $\mathcal F \varphi \in PC_{\omega}.$ 

Next, the proof is completed by showing  that $\mathcal F$ is a contraction mapping. 
If $\varphi^1, \varphi^2\in PC_{\omega},$ then
\begin{eqnarray*}
||\mathcal F\varphi^1(t)-\mathcal F\varphi^2(t)||&=&\sum^{m}_{i=1}\left|(\mathcal F\varphi^1)_i(t)-(\mathcal F\varphi^2)_i(t)\right|\\&\leq&\sum^{m}_{i=1}\left\{\int_0^{\omega}|\mathcal G_{i}(t,s)|
\left[\displaystyle\sum_{j=1}^{m}L_{j}^{f}|b_{ij}||\varphi_{j}^1(s)-\varphi_{j}^2(s)|\right. \right.\\  &&+ \left. \lambda
\sum_{j=1}^{m} L_{j}^{g}|c_{ij}||\varphi_{j}^1(s)-\varphi_{j}^2(s)|\right]ds  \\ &&+ \left. \ell\sum_{k=1}^{p}|\mathcal G_i(t,\tau_{k})||\varphi_{i}^1(\tau_{k}^{-})-\varphi_{i}^2(\tau_{k}^{-})|\right\}\\&\leq&\mathcal R\sum^{m}_{i=1}\left\{\int_0^{\omega}\left[\sum_{j=1}^{m}L_{i}^{f}|b_{ji}||\varphi_{i}^1(s)-\varphi_{i}^2(s)|\right.\right.\\&&\left.\left.+\sum_{j=1}^{m} L_{i}^{g}|c_{ji}||\varphi_{i}^1(\beta(s))-\varphi_{i}^2(\beta(s))|\right]ds+\ell\sum_{k=1}^{p}|\varphi_{i}^1(\tau_{k}^{-})-\varphi_{i}^2(\tau_{k}^{-})|\right\}\\&\leq& \mathcal R\left(\int_0^{\omega}\left[\mu||\varphi^1(s)-\varphi^2(s)||+k_2||\varphi^1(\beta(s))-\varphi^2(\beta(s))||\right]ds\right.\\&&\left.+\ell\sum_{k=1}^{p}||\varphi^1(\tau_k^{-})-\varphi^2(\tau_k^{-})||\right).
\end{eqnarray*}
Hence,
\begin{eqnarray*}
||\mathcal F\varphi^1-\mathcal F\varphi^2||_0&\leq&\mathcal R\Big(\omega\left(\mu+\lambda k_2\right)+\ell p\Big)||\varphi^1-\varphi^2||_0.
\end{eqnarray*}
It follows from the  condition (H7) that,  $\mathcal F$ is a contraction mapping in $PC_{\omega}.$ Consequently, by using Banach fixed point theorem, $\mathcal F$ has a unique fixed point $\varphi^*\in PC_{\omega},$
such that $\mathcal F\varphi^*=\varphi^*.$ This completes the proof. $\Box$

\begin{theorem}\label{pgs}
Assume that conditions $\rm{(H1)-(H7)}$ are valid. 
 Then the periodic solution of (\ref{ene}) is globally asymptotically stable.
\end{theorem}

\noindent {\bf Proof.}
By Theorem \ref{ups}, we know that (\ref{ene}) has an $\omega-$periodic solution $z^*(t)=(z_{1}^*,\cdots,z_{m}^*)^T.$ Suppose that $z(t)=(z_1,\cdots,z_m)^T$ is an arbitrary solution of (\ref{ene}) and let $z(t)=z(t)-z^*(t)=(z_1-z_{1}^*,\cdots,z_m-z_{m}^*)^T.$ Then, similar to the proof of Theorem \ref{tgs}, one can show that it is  globally asymptotically stable.

\section{ Numerical simulations}\label{simulations}
In this section, we give  examples with numerical simulations to illustrate the theoretical  results of the paper. In what follows, let  $\theta_{k}=k,\ \tau_{k}=(\theta_{k}+\theta_{k+1})/2=(2k+1)/2,\ k\in \mathbb N$ be the sequence of the change of constancy for the argument and the sequence of impulsive action, respectively. 

Consider the following recurrent neural networks:  
{\footnotesize\begin{eqnarray}\label{ex1}
\left\{\begin{array}{lll}
\frac{dx(t)}{dt}&=&-\left(
\begin{array}{ccc} 5\times10^{-1} & 0 \\ 0 & 5\times10^{-1} 
\end{array}\right)\left(
\begin{array}{ccc} x_{1}(t)  \\ x_{2}(t) 
\end{array}\right)+ \left(
\begin{array}{ccc} 10^{-4} & \ 2\times10^{-4} \\ 10^{-4} & 3\times10^{-4}
\end{array}\right)\left(
\begin{array}{ccc} \tanh(\frac{x_{1}(t)}{10})  \\ \tanh(\frac{3x_{2}(t)}{10})  
\end{array}\right) \\&&+
\left(
\begin{array}{ccc} 2\times 10^{-2} & 3\times 10^{-3} \\ 3\times 10^{-3} & 5\times 10^{-3} 
\end{array}\right)\left(
\begin{array}{ccc} \tanh(\frac{x_{1}(\beta (t))}{5})  \\ \tanh(\frac{x_{2}(\beta (t))}{5})  
\end{array}\right)+\left(
\begin{array}{ccc} 1  \\ 1 
\end{array}\right), \ t\neq \tau_k \\ 
\Delta x(t)&=& \left(\begin{array}{ccc} I_k(x_{1}(\tau_{k}^-))\\ I_k(x_{2}(\tau_{k}^-)) 
\end{array}\right)= \left(\begin{array}{ccc} \frac{x_{1}(\tau_{k}^-)}{40}+\frac{1}{2}\\ \frac{x_{2}(\tau_{k}^-)}{40}+\frac{1}{2} 
\end{array}\right),\ t=\tau_k, \quad k=1,2,\cdots,
\end{array}\right.
\end{eqnarray}}

By simple calculation, one can see that the corresponding parameters in the conditions of Theorems \ref{eus}, \ref{tgs}, \ref{ups}, \ref{pgs} are $k_{1}=0.5001,\ k_{2}=0.0046,\ L_{1}^{f}=0.1,\ L_{2}^{f}=0.3,\ L_{1}^{g}=L_{2}^{g}=0.2,\ \ell=0.0250,\  \overline{\theta}=\underline{\tau}=1,\ p=p_{1}=1,\ \gamma=0.5,\ \lambda=9.6421,\ \mu=0.00015,\ \omega=1,\ \mathcal R=2.5415, \ \alpha_{1}=0.1766.$ For these values, we can check that $(H3)=0.9032 < 1,\  (H4)=0.8963 < 1, (H5)=0.4308 > 0$ and $\alpha_{1}=0.1766 < 1.$ So, it is easy to verify that (\ref{ex1}) satisfies the conditions of these theorems. Hence, the system   of  (\ref{ex1}) has a  1-periodic solution which is globally asymptotically stable. Specifically, the simulation results with some  initial points  are shown in Fig. 1. and Fig. 2. We deduce that  the non-smoothness at $\theta_k, \ k\in \mathbb N$ is not  seen by numerical simulations due to the choosing the parameters small enough to satisfy the theorems. Hence, the {\it smallness} hides the {\it non-smoothness}. 

\begin{figure}[!htp] 
\centering
\includegraphics[width=12cm,height=9cm]{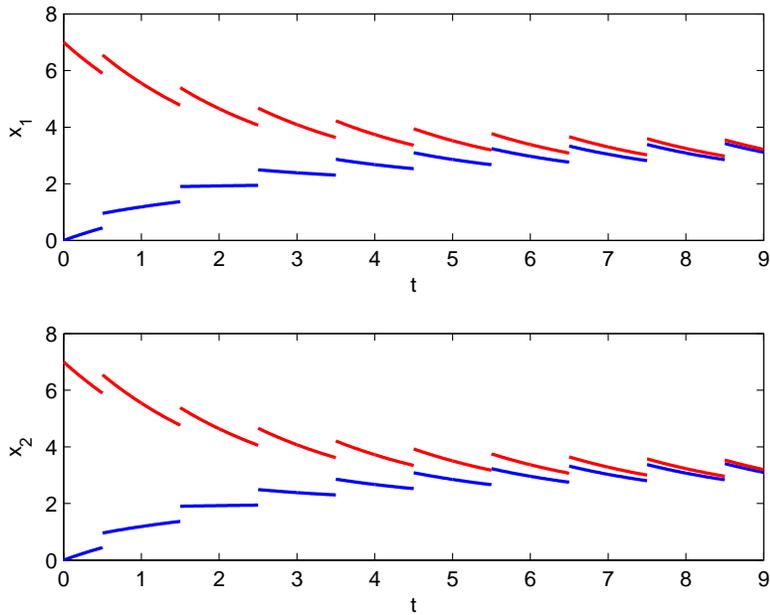}
\caption{ Transient behavior of the recurrent neural networks for the system (\ref{ex1}) with the initial points $[0, 0]^{T}$ and $[7, 7]^{T}.$}
\label{fig1}
\end{figure}

\begin{figure}[!htp] 
\centering
\includegraphics[width=12cm,height=9cm]{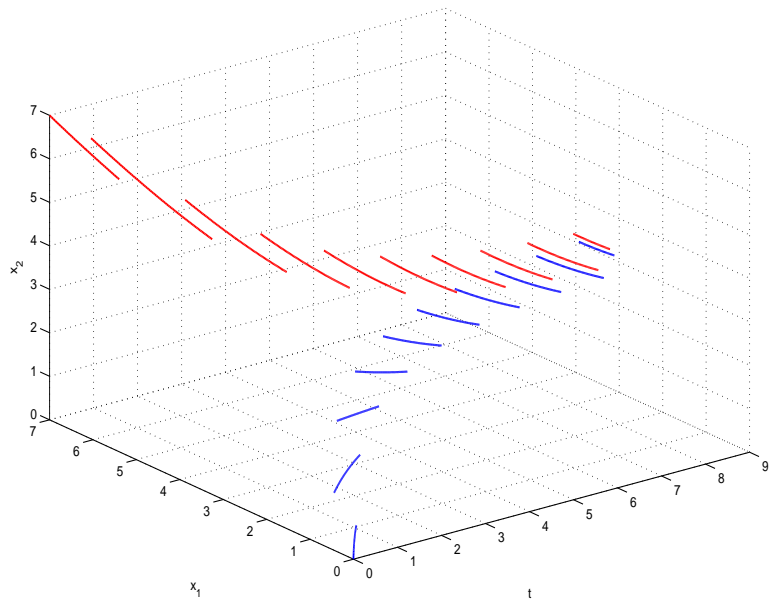}
\caption{ Eventually  1-periodic solutions of system (\ref{ex1}) with the initial points $[0, 0]^{T}$ and $[7, 7]^{T}.$
}\label{fig2}
\end{figure}

\newpage
On the other hand, in the following example, we illustrate a globally stable equilibrium appearance for our system of differential equations:

{\footnotesize\begin{eqnarray}\label{ex3}
\left\{\begin{array}{lll}
\frac{dx(t)}{dt}&=&-\left(
\begin{array}{ccc} 5\times10^{-1} & 0 \\ 0 & 5\times10^{-1} 
\end{array}\right)\left(
\begin{array}{ccc} x_{1}(t)  \\ x_{2}(t) 
\end{array}\right)+ \left(
\begin{array}{ccc} 10^{-4} & \ 2\times10^{-3} \\ 10^{-4} & 3\times10^{-3}
\end{array}\right)\left(
\begin{array}{ccc} \tanh(\frac{x_{1}(t)}{10})  \\ \tanh(\frac{3x_{2}(t)}{10})  
\end{array}\right) \\&&+
\left(
\begin{array}{ccc} 2\times 10^{-2} & 3\times 10^{-2} \\ 3\times 10^{-2} & 5\times 10^{-2} 
\end{array}\right)\left(
\begin{array}{ccc} \tanh(\frac{x_{1}(\beta (t))}{5})  \\ \tanh(\frac{x_{2}(\beta (t))}{5})  
\end{array}\right)+\left(
\begin{array}{ccc} 1  \\ 1 
\end{array}\right), \ t\neq \tau_k \\ 
\Delta x(t)&=& \left(\begin{array}{ccc} I(x_{1}(\tau_{k}^-))\\ I(x_{2}(\tau_{k}^-)) 
\end{array}\right)= \left(\begin{array}{ccc} \frac{(x_{1}(\tau_{k}^-)-x_{1}^*)^{2}}{30}\\ \frac{(x_{2}(\tau_{k}^-)-x_{1}^*)^{2}}{30} 
\end{array}\right),\ t=\tau_k, \quad k=1,2,\cdots,
\end{array}\right.
\end{eqnarray}}
where $x_{1}^*=2.0987,\ x_{2}^*=2.1577.$ One can check that  the point  $x^*=(x_{1}^*,x_{2}^*)$ satisfies the algebraic system
\begin{eqnarray}\label{alsy}
-a_{i}x_{i}^*+\displaystyle\sum_{j=1}^{2}b_{ij}f_j(x_{j}^*) + \sum_{j=1}^{2}c_{ij}g_j(x_{j}^*))+d_{i}=0,
\end{eqnarray}
approximately. And it is clear that   $I(x_{i}^*)=0$ for  $i=\overline{1,2}.$ By simple calculation, we can see that all conditions of  Theorem \ref{eueq} are satisfied and the point $x^{*}$ is a solution of (\ref{alsy}), approximately  with the error, which is less than $10^{-11}$(evaluated by MATLAB). 

The simulation, where the initial value is chosen as $[10, 10]^{T}$ , is shown in Fig. 3 and it illustrates that all trajectories  converge to  $x^{*}.$

\begin{figure}[!htp] 
\centering
\includegraphics[width=12cm,height=7cm]{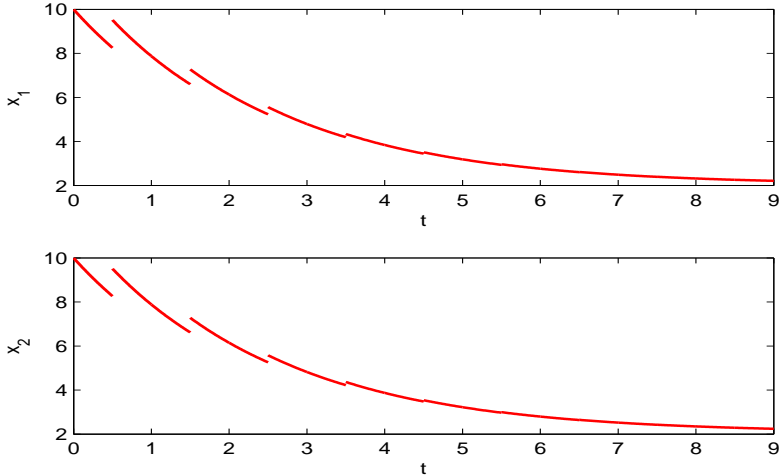}
\caption{The first and the second coordinates of the solution for the the system (\ref{ex3}) with the initial point $[10, 10]^{T}$ approaches $x_{1}^*$ and $x_{2}^*$, respectively, as time increases.} 
\label{fig5}
\end{figure}

Now, let us take the parameters so that the non-smoothness  can also be seen. Consider the following recurrent neural networks with non-smooth and impact activations:

\begin{eqnarray}\label{ex2}
\left\{\begin{array}{lll}
\frac{dx(t)}{dt}&=&-\left(
\begin{array}{ccc} 20 & 0 \\ 0 & 8 
\end{array}\right)\left(
\begin{array}{ccc} x_{1}(t)  \\ x_{2}(t) 
\end{array}\right)+ \left(
\begin{array}{ccc} 10 & \ 3 \\ 5 & 1
\end{array}\right)\left(
\begin{array}{ccc} \tanh(x_{1}(t))  \\ \tanh(x_{2}(t))  
\end{array}\right) \\&&+
\left(
\begin{array}{ccc} 20 & 1 \\ 8 & 1 
\end{array}\right)\left(
\begin{array}{ccc} \tanh(x_{1}(\beta (t)))  \\ \tanh(x_{2}(\beta (t)))  
\end{array}\right)+\left(
\begin{array}{ccc} 1  \\ 1 
\end{array}\right), \ t\neq \tau_k \\ 
\Delta x(t)&=& \left(\begin{array}{ccc} I_k(x_{1}(\tau_{k}^-))\\ I_k(x_{2}(\tau_{k}^-)) 
\end{array}\right)= \left(\begin{array}{ccc} \frac{x_{1}(\tau_{k}^-)}{3}+\frac{1}{6}\\ \frac{x_{2}(\tau_{k}^-)}{3}+\frac{1}{6} 
\end{array}\right),\ t=\tau_k, \quad k=1,2,\cdots,
\end{array}\right.
\end{eqnarray}

Clearly, one can see that our parameters  are big now. Therefore,  the system of equations (\ref{ex2}) does not satisfy the conditions of the theorems. However, we can see the {\it non-smoothness} of the solution with the initial value $[0, 0]^{T},$ which is illustrated by simulations in Fig.3. and Fig.4.

\begin{figure}[ht] 
\centering
\includegraphics[width=12cm,height=9cm]{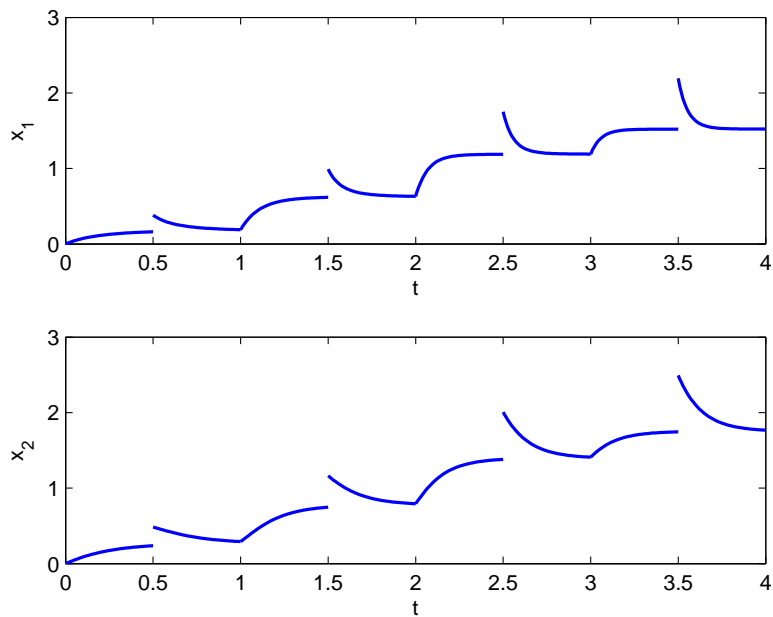}
\caption{ The {\it impact} and {\it non-smoothness} are seen at discontinuity points $\tau_{k} \ (0.5; 1.5; 2.5; 3.5)$ and at  switching points $\theta_{k} \ (1; 2; 3),$ respectively.}
\label{fig3}
\end{figure}

\begin{figure}[ht] 
\centering
\includegraphics[width=12cm,height=9cm]{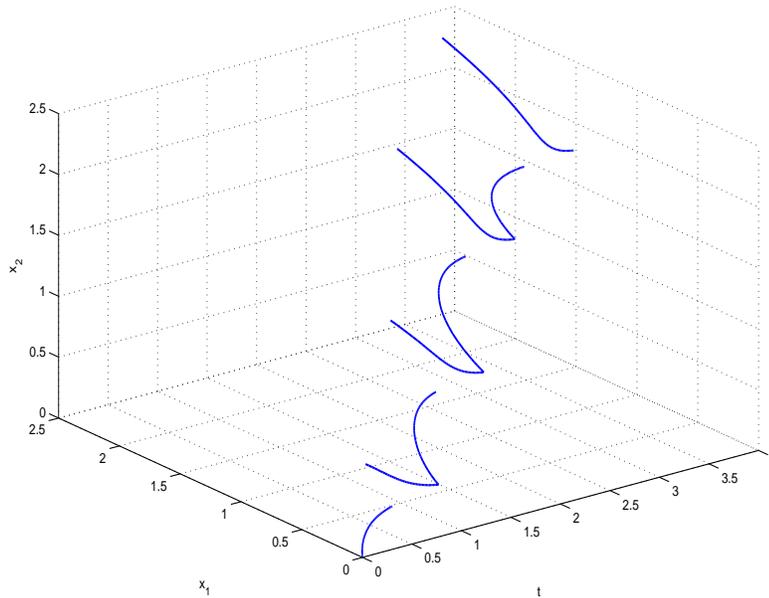}
\caption{ Eventually  1-periodic solutions of system (\ref{ex2})  with the initial point $[0, 0]^{T}.$
}\label{fig4}
\end{figure}

\section{Conclusions}

This  is the first time that global asymptotic stability of periodic solutions for recurrent neural networks with  both impulses and   piecewise constant delay is considered.  Furthermore, our model gives new ideas not only from the implementation point of view, but also from the system of differential equations. In other words, we develop differential equations with piecewise constant argument  to  a new class of system, so called impulsive differential equations with piecewise constant delay. For applications, we have also nice properties   on the system of equations that   the moments of discontinuity $\tau_k$ and switching moments of constancy of arguments $\theta_k$  are not   related to each other. That is,  our investigations are more  applicable to the  real world problems like recurrent neural networks.  Finally, the  results given in this paper could be developed for  more complex systems \cite{ak8}.

%{\bfseries{Acknowledgement}}
%
%The authors wish to express their sincere gratitude to the referees for the useful suggestions, which helped them to improve the paper.

\end{document}